\documentclass[final,1p,times]{elsarticle}
\usepackage{graphicx,subfigure,epsfig}
\usepackage{amsmath,mathrsfs,amssymb,amsthm,dsfont,bm,amsfonts}         %
\usepackage{color}
\usepackage{amsmath}
\allowdisplaybreaks[4]

\newtheorem{theorem}{Theorem}[section]
\newtheorem{re}{Remark}[section]

\newtheorem{lemma}{Lemma}[section]

\newtheorem{assumption}{Assumption}[section]
\newtheorem{alg}{Algorithm}

\newcommand{\p}{\partial}
\newcommand{\f}{\frac}

\newcommand{\up}{\textup}
\newcommand{\tr}{\triangle}

\newcommand{\mbf}{\mathbf}

\journal{*******************}
\begin{document}

\begin{frontmatter}

\title{Hybrid mixed discontinuous Galerkin finite element method for incompressible miscible displacement problem}

\author[a]{Jiansong Zhang\corref{zhang}}
\ead{jszhang@upc.edu.cn}

\author[a]{Yun Yu}
\ead{yuyun19970321@163.com}

\author[b]{Jiang Zhu}
\ead{jiang@lncc.br}

\author[a]{Rong Qin}
\ead{qr-920diana@126.com}

\author[a]{Yue Yu}
\ead{m18766215811@163.com}

\author[c]{Maosheng Jiang}
\ead{msjiang@qdu.edu.cn}

\address[a]{Department of Applied Mathematics,
        China University of Petroleum, Qingdao 266580, China.}
\address[b]{Laborat\'orio Nacional de Computa\c{c}\~ao Cient\'{\i}fica,   Petr\'opolis 25651-075, RJ, Brazil}
\address[c]{School of Mathematics,
       Qingdao University, China.}
\cortext[zhang]{Corresponding author.}

\begin{abstract}
A new hybrid mixed discontinuous Galerkin finite element (HMDGFE) method is constructed for incompressible miscible displacement problem. In this method, the hybrid mixed finite element (HMFE) procedure is considered to solve pressure and velocity equations, and a new hybrid mixed discontinuous Galerkin procedure is constructed to solve the concentration equation with upwind technique. Compared with other traditional discontinuous Galerkin methods, the new method can reach global systems with less unknowns and sparser stencils. The consistency and conservation of the method are analyzed, the stability and optimal error estimates are also derived by the new technique. 
\end{abstract}


\end{frontmatter}

\section{Introduction}
 In this paper, we consider a new numerical method for simulating two-phase (water and oil) miscible displacement problem in porous media  (see \cite{model1,model2,model3,model3,model4,model5}):
\begin{equation}\label{eq1}
\begin{aligned}
&\nabla\cdot \mathbf{u} = q,\quad \mathbf{u}=\frac{-k^{\ast}}{\mu(c)}\nabla p,\quad x\in \Omega,\quad 0 \le t \leq T, \\
&\phi\frac{\partial c}{\partial t}+\nabla\cdot(\mathbf{u}c-D(\mathbf{u})\nabla c)=qc^{\ast},\quad x\in \Omega,\quad 0 \le t \leq T,
\end{aligned}
\end{equation}
with the following initial-boundary conditions
\begin{equation}\label{eq2}
\begin{aligned}
\mathbf{u}\cdot \mathbf{n}&=0,\quad x\in \partial\Omega,\quad 0 \le t \le T,\\
D(\mathbf{u})\cdot \mathbf{n}&=0,\quad x\in \partial\Omega,\quad 0 \le t \le T,\\
c(x,0)&=c_{0}(x),\quad x\in \Omega,
\end{aligned}
\end{equation}
where $\Omega$ is a bounded polygonal or polyhedral domain in $R^d(d=2,3)$;  $\mathbf{n}$ denotes the unit outard normal vector  of boundary $\partial\Omega$; the functions $p$ and $\mathbf{u}$ denote the pressure and Darcy velocity; $c$ is the concentration of the fluid mixture; $c^{\ast}$ is the concentration in the external flow, which must be specified at injection points
$(q > 0)$ and is assumed to be equal to $c$ at production points $(q < 0)$; $k^{\ast}$ and $\phi$ are the permeability and
the porosity of the rock, respectively; $\mu(c)$ is the viscosity of the fluid.
Diffusion coefficient $D(\mathbf{u})=\phi[d_mI+|\mathbf{u}|(d_IE(\mathbf{u})+d_tE^{\perp}(\mathbf{u}))]$ comes from two aspects: small molecule diffusion of oil field scale problem, and speed-related diffusion in petroleum engineering, here the $d \times d$ matrix $E(\mathbf{u})=(u_iu_j/|\mathbf{u}|^2)_{d \times d}$  represents orthogonal projection along the velocity vector, and $E^{\perp}(\mathbf{u}) = I-E(\mathbf{u})$, $d_m$, $d_l$ and $d_t$ are the molecular diffusion, longitudinal and transverse dispersion coefficients, respectively.

We first consider numerical method for the pressure and velocity equations. 
As we know, classical mixed finite element methods \cite{model1,model2,mixed3,mixed4} can achieve more accurate approximation of the velocity. However, these methods have some disadvantages: it can cause saddle point problems, the finite element spaces require the LBB condition and so on. To avoid these problems, hybrid mixed element methods \cite{hmixed1,hmixed2,hmixed3} are widely used, because they can eliminate the continuity of the normal component of the velocity over each element interface and lead to a symmetric positive definite system. So we consider the hybrid mixed element method to solve Darcy velocity $\mathbf{u}$ and pressure $p$. An important advantage of this method is that subproblems for velocity and pressure can be solved at element level and these variables are eliminated in favor of the Lagrange multiplier, which is identified as pressure trace at the element interfaces, and the global system involves only the degrees of freedom associated with the multiplier, significantly reducing the computational cost.

Then, we turn to the numerical approximation for the concentration $c$. Generally speaking, the standard Galerkin finite element method does not solve the the concentration problems well, especially for the convection-dominated and  discontinuous cases. Discontinuous Galerkin (DG) finite element methods \cite{ar82}  have many advantages. For example, test functions across the finite element interfaces have no continuity requirement, thus the spaces are easy to construct, and the use of highly nonuniform and unstructured meshes is permitted. The combined traditional mixed element methods with DG procedure for miscible displacement problem were presented in \cite{hgmixed1,hgmixed2,hgmixed3,hgmixed4,hgmixed5}. Moreover, Zhu, Zhang etc. considered DG methods to other different problems, see \cite{zhu3,zhu4,zhu1,zhu2}. Howerver, the traditional DG methods have some disadvantages: the overall number of unknowns is increased substantially compared with a standard conforming discretization, and the resulted linear systems are much less sparse. In order to overcome these problems, the Lagrange multiplier in the DG methods is introduced, which couples hybrid mixed methods and DG methods naturally, and gets a stable mixed hybrid DG procedure.  In this way, one can reach global systems with less unknowns than other DG methods.

In this article, we will combine the hybrid mixed finite element method with the discontinuous Galerkin method for solving miscible displacement problem in porous media. Here we will consider the hybrid mixed finite element procedure  for pressure and velocity, and then, based on the combination of a hybrid mixed element method and a discontinuous Galerkin method, we construct a new hybrid mixed discontinuous Galerkin finite element method for the concentration with upwind technique. Compared with the existing combined methods as in \cite{hgmixed1,hgmixed2,hgmixed3,hgmixed4,hgmixed5}, the main contribution of this article is to construct a new hybrid mixed DG method for the concentration and to present the convergence analysis by a different technique from the ones of the traditional DG methods. Meantime, the optimal error estimate is derived  under the case that the diffusion coefficient includes the molecular diffusion and dispersion. 

For the convenience  of analysis,  we  make the following assumptions on the coefficient parameters and the regularities of the solution of the system \eqref{eq1}.
\begin{assumption}\label{eq3}
Assume that the following parameters $\phi(x)$, $\frac{k(x)^{\ast}}{\mu(c)}$ and $q(x,t)$ are bounded as follows:
\begin{equation}
\begin{aligned}
0<\phi_{\ast}\leq\phi(x)\leq\phi^{\ast},\quad
0<a_{\ast}\leq\frac{k^{\ast}}{\mu(c)}\leq a^{\ast},\quad | q|\leq C,
\end{aligned}
\end{equation}
where  $\phi_{\ast}$, $\phi^{\ast}$, $a_{\ast}$, $a^{\ast}$ and $C$ is a positive constant.
\end{assumption}
\begin{assumption}\label{eq4}
Assume that the solution ($p$, $\mathbf{u}$, $c$) of the system \eqref{eq1} satisfies the following regularities:
\begin{equation}
\begin{aligned}
&p\in L^{\infty}(0,T;H^{k+1}(\Omega)),\quad \mathbf{u}\in L^{\infty}(0,T;H^{k+1}(\Omega))\cap L^{\infty}(0,T;L^{\infty}(\Omega)),\\
& \bm\sigma\in L^{\infty}(0,T;H^{k+1}(\Omega))\cap L^{\infty}(0,T;L^{\infty}(\Omega)),\quad c\in L^{\infty}(0,T;H^{k+1}(\Omega))\cap H^{2}(0,T;L^{2}(\Omega)).
\end{aligned}
\end{equation}
\end{assumption}

\section{The formulation of the method}

Let $\mathcal{T}_h$ be a quasi-uniform regular partition of $\Omega$ with $\mathcal{T}_h = \{K_1, K_2, . . . , K_N \}$, and denote $\partial\mathcal{T}_h=\cup_{K\in\mathcal{T}_h}\{e|e\in\partial K\}$ to be the set of all cell edges. The velocity vector field $\mathbf{u}$ induces a natural splitting of element boundaries into inflow and outflow parts, i.e., we denote the outflow boundary $\partial K^{out}=\{x\in\partial K:\mathbf{u}\cdot\mathbf{n}>0\}$ and $\partial K^{in}=\partial K\backslash\partial K^{out}$, where $\mathbf{n}$ denotes the unit normal direction of $\p K$. The unions of the element inflow and outflow boundaries are $\partial\mathcal{T}_h^{in}$ and $\partial\mathcal{T}_h^{out}$, respectively. And $\partial\Omega^{in}$ and $\partial\Omega^{out}$ are the inflow and outflow regions of the boundary $\partial\Omega$. Furthermore, let $h_e = \textrm{diam}(e)$ and $h=\max_e(h_e)$ for all $e\in\partial\mathcal{T}_h$.

In this paper, we need to introduce the piecewise Sobolev spaces
\begin{equation}\nonumber
H^s(\mathcal{T}_h)=\{v\in L^2(\Omega):v|_K\in H^s(K),K\in\mathcal{T}_h\},\quad s\ge 0,
\end{equation}
and
\begin{equation}\nonumber
\begin{aligned}
&L^2(\partial\mathcal{T}_h)=\{v\in L^2(e),\forall e\in\partial\mathcal{T}_h\}.
\end{aligned}
\end{equation}
Define some inner products as follows:
\begin{equation}\nonumber
\begin{aligned}
&(u,v)_{K}=\int_{K}uvdx\quad (u,v)_{\mathcal{T}_h}=\sum\limits_{K\in\mathcal{T}_h}(u,v)_K,\\
 & \langle u,v\rangle_{\p  K}=\int_{\p K}uvds\quad \langle u,v\rangle_{\partial\mathcal{T}_h}=\sum\limits_{K\in\mathcal{T}_h}\langle u,v\rangle_{\p K},
\end{aligned}
\end{equation}
and the norms $\|\cdot\|_{\mathcal{T}_{h}}=\sqrt{(\cdot,\cdot)_{\mathcal{T}_h}}$ and  $|\cdot|_{\partial\mathcal{T}_{h}}=\sqrt{\langle\cdot,\cdot\rangle_{\partial\mathcal{T}_h}}$.

Introduce the discrete approximate spaces  $\Psi_h$, $\Lambda_h$, $\Pi_h$ and $\Sigma_h$
\begin{equation}\nonumber
\begin{aligned}
&\Psi_h=\{v\in H^k(\mathcal{T}_h):v|_K\in P_k(K),K\in\mathcal{T}_h\},\\
&\Lambda_h=\{v\in L^2(\mathcal{T}_h):v|_K\in P_k(K),K\in\mathcal{T}_h\},\\
&\Theta_h=\{v\in [L^2(\mathcal{T}_h)]^d:v|_K\in RT_k(K),K\in\mathcal{T}_h\},\\
&\Sigma_h=\{v\in L^2(\partial\mathcal{T}_h):v|_e\in P_k(e),e\in\partial\mathcal{T}_h\},
\end{aligned}
\end{equation}
where $P_k(K)$, $P_k(e)$ are the spaces of polynomial functions of degree at most $k$ for each $K\in\mathcal{T}_h$ and each $e\in\partial\mathcal{T}_h$, respectively, and $RT_k(K)=[P_k(K)]^d\oplus xP_k(K)$ denotes the Raviart-Thomas element space as in \cite{hmixed1,hmixed2,hmixed3}.  

Let $M>0$ be a positive integer and $\Delta t=T/M$ be the time size, and  denote $t^{n}=n\Delta t$, $n=0,1,\cdots, M$.  Next, we will formulate our method for miscible displacement problem
in porous media.

\subsection{Revision  of hybrid mixed finite element method for pressure and velocity}

In this subsection, we give the hybrid mixed finite element  method for pressure and velocity as  in \cite{hgmixed2}. Here we use completely discontinuous piecewise polynomial functions and ensure the continuity of the normal fluxes over internal interfaces by adding Lagrangian multiplier. The hybrid mixed finite element formulation can be written as below:
\begin{alg}\label{alg1}
For given approximate value  $c_{h}^{n-1}$, seek $(\mathbf{u}_h^{n}, p_h^{n}, \lambda_h^{n}) \in  \Theta_h \times \Lambda_h \times \Sigma_h$ such that
\begin{equation}\label{eq5}
\begin{aligned}
(a(c_{h}^{n-1})\mathbf{u}_h^{n},\bm\tau_h)_{\mathcal{T}_h}-(p_h^{n},\nabla\cdot\bm\tau_h)_{\mathcal{T}_h}+\langle\lambda_h^{n},\bm\tau_h\cdot \mbf{n}\rangle_{\partial\mathcal{T}_h}&=0,\quad&\forall \bm\tau_h \in \Theta_h,\\
(\nabla\cdot\mathbf{u}_h^{n},v_h)_{\mathcal{T}_h}&=(q^{n},v_h)_{\mathcal{T}_h},\quad&\forall v_h \in \Lambda_h,\\
\langle\mathbf{u}_h^{n}\cdot \mbf{n},\mu_h\rangle_{\partial\mathcal{T}_h}&=0,\quad&\forall \mu_h \in \Sigma_h,
\end{aligned}
\end{equation}
where $a(c_{h}^{n-1})=\mu(c_{h}^{n-1})/k^{\ast}$.
\end{alg}

Define the bilinear form:
\[
\begin{aligned}
&B_{\mathbf{u}}(c_h^{n-1};(\mathbf{u}_h^{n},p_h^{n},\lambda_h^{n}),(\bm\tau_h,v_h,\mu_h))\\&\quad:=(a(c_h^{n-1})\mathbf{u}_h^{n},\bm\tau_h)_{\mathcal{T}_h}+(\nabla p_h^{n},\bm\tau_h)_{\mathcal{T}_h}+(\mathbf{u}_h^{n},\nabla v_h)_{\mathcal{T}_h}\\
&\quad\quad+\langle\lambda_h^{n}-p_h^{n},\bm\tau_h\cdot \mbf{n}\rangle_{\partial\mathcal{T}_h}+\langle\mathbf{u}_h^{n}\cdot \mbf{n},\mu_h-v_h\rangle_{\partial\mathcal{T}_h}.
\end{aligned}
\]
We can rewrite \eqref{eq5} into the following equivalent form:
\begin{alg}[HMFE Algorithm]\label{alg2}
 For given $c_h^{n-1}$, find $(\mathbf{u}_h^{n}, p_h^{n}, \lambda_h^{n}) \in \Theta_h \times \Lambda_h \times \Sigma_h$ such that
\begin{equation}\label{eq6}
\begin{aligned}
B_{\mathbf{u}}(c_h^{n-1};(\mathbf{u}_h^{n},p_h^{n},\lambda_h^{n}),(\bm\tau_h,v_h,\mu_h))=-(q^{n},v_h)_{\mathcal{T}_h}, \quad  \forall (\bm\tau_h, v_h, \mu_h) \in \Theta_h \times \Lambda_h \times \Sigma_h.
\end{aligned}
\end{equation}
\end{alg}
\subsection{Hybrid mixed discontinuous Galerkin finite element scheme for the concentration}
Set $\bm\sigma=-D(\mathbf{u})\nabla c$. We can rewrite the concentration equation of \eqref{eq1} into the following first order partial differential equations:
\begin{equation}\nonumber
\begin{aligned}
\bm\sigma+D(\mathbf{u})\nabla c&=0,\\
\phi\frac{\partial c}{\partial t}+\nabla\cdot(\bm\sigma+\mathbf{u}c)&=qc^{\ast}.
\end{aligned}
\end{equation}
Multiplying the above two equations by the test functions $\bm\tau$ and $v$ respectively,  and adding an upwind stabilization term, we can reach the mixed variational formulation for the concentration equation.
\begin{equation}\label{eq7}
\begin{aligned}
&(D(\mathbf{u})^{-1}\bm\sigma,\bm\tau)_{\mathcal{T}_h}+(\nabla c,\bm\tau)_{\mathcal{T}_h}=0,\\
&(\phi\frac{\partial c}{\partial t},v)_{\mathcal{T}_h}+(\nabla\cdot(\bm\sigma+\mathbf{u}c),v)_{\mathcal{T}_h}+\langle\mathbf{u}\cdot\mathbf{n}(c^{+}-c),v\rangle_{\partial\mathcal{T}^{in}_h}=(qc^{\ast},v)_{\mathcal{T}_h},
\end{aligned}
\end{equation}
where $c^{+}:=c|_{\partial K^{+}}$ denotes the upwind value at the upwind element $K^+$, that is, the element attached to $e$ where $\mathbf{u}\cdot\mathbf{n}=\mathbf{u}\cdot\mathbf{n}_{e}\geq 0$. To incorporate the boundary condition, we define $c^+=0$
on $\partial\Omega^{in}$. After integration by parts, we can reach
\begin{equation}\nonumber
\begin{aligned}
&(D(\mathbf{u})^{-1}\bm\sigma,\bm\tau)_{\mathcal{T}_h}+(\nabla c,\bm\tau)_{\mathcal{T}_h}=0,\\
& (\phi\frac{\partial c}{\partial t},v)_{\mathcal{T}_h}-(\bm\sigma+\mathbf{u}c,\nabla v)_{\mathcal{T}_h}+\langle\mathbf{u}\cdot\mathbf{n}c^{+},v\rangle_{\partial\mathcal{T}_h^{in}}+\langle\mathbf{u}\cdot\mathbf{n}c,v\rangle_{\partial\mathcal{T}_h^{out}}=(qc^{\ast},v)_{\mathcal{T}_h}.
\end{aligned}
\end{equation}

Introduce the upwind value as a new variable $\lambda^{c}:=c^+$, and define $\hat{c}$ as follows: for any $K\in\mathcal{T}_{h}$
\begin{equation}\label{flux}
\hat{c}:=\left\{
\begin{aligned}
&\lambda^c,& \quad e\subset\partial K^{in},\\
&c,& \quad e\subset\partial K^{out}.\\
\end{aligned}
\right.
\end{equation}
Using the fact that $\lambda^{c}=\hat{c}=c^+$ on both sides of $e$, we can give the corresponding fully discrete hybrid mixed discontinuous Galerkin method for the concentration:
 \begin{alg}\label{alg3}
 For given approximate value $\mathbf{u}_{h}^{n}$, seek $(\bm\sigma_h^{n}, c_h^{n}, \lambda_h^{c,n}) \in  \Theta_h \times \Lambda_h \times \Sigma_h$
\begin{equation}\label{con_pro_3}
\begin{aligned}
&(D(\mathbf{u}^{n}_h)^{-1}\bm\sigma^{n},\bm\tau_h)_{\mathcal{T}_h}-(c_h^{n},\nabla\cdot\bm\tau_h)_{\mathcal{T}_h}+\langle\lambda^{c,n}_h,\bm\tau_h\cdot\mathbf{n}\rangle_{\partial\mathcal{T}_h}=0,\quad\bm\tau_h\in\Theta_h,\\
\end{aligned}
\end{equation}
\begin{equation}\label{con_pro_4}
\begin{aligned}
(\phi\bar{\partial}c_h^{n},v_h)_{\mathcal{T}_h}&+(\nabla\cdot(\mathbf{u}^{n}_hc_h^{n})+\nabla\cdot\bm\sigma_h^{n},v_h)_{\mathcal{T}_h}
\\
&+\langle\mathbf{u}_h^{n}\cdot\mathbf{n}(\lambda_h^{c,n}-c_h^{n}),v_h\rangle_{\partial\mathcal{T}_h^{in}}=(q^{n}c^{*,n}_{h},v_h)_{\mathcal{T}_h},\quad v_h\in\Lambda_h,\\
\end{aligned}
\end{equation}
\begin{equation}\label{con_pro_1}
\langle\bm\sigma^{n}_h\cdot\mathbf{n},\mu_h\rangle_{\partial\mathcal{T}_h}=0,\quad\mu_h\in\Sigma_h,
\end{equation}
\begin{equation}\label{con_pro_2}
\langle\mathbf{u}^{n}_h\cdot\mathbf{n}\hat{c}_h^{n},\mu_h\rangle_{\partial\mathcal{T}_h}=0,\quad\mu_h\in\Sigma_h,
\end{equation}
where $\bar{\partial}c_h^{n}=(c_h^{n}-c_h^{n-1})/\Delta t$.
\end{alg}

Define the bilinear form:
\[
\begin{aligned}
&B_c(\mathbf{u}_h^{n};(\bm\sigma_h^{n},c_h^{n},\lambda_h^{c,n}),(\bm\tau_h,v_h,\mu_h))
\\&\quad:=-\frac{1}{\Delta t}(\phi c_h^{n},v_h)_{\mathcal{T}_h}+(\bm\sigma_h^{n}+\mathbf{u}^{n}_{h}c_h^{n},\nabla v_h)_{\mathcal{T}_h}\\&\quad\quad+(\frac{1}{D(\mathbf{u}_h^{n})}\bm\sigma_h^{n}+\nabla c_h^{n},\bm\tau_h)_{\mathcal{T}_h}
+\langle\lambda_h^{c,n}-c_h^{n},\bm\tau_h\cdot\mathbf{n}\rangle_{\partial\mathcal{T}_h}\\
&\quad\quad+\langle\bm\sigma_h^{n}\cdot\mathbf{n}+\mathbf{u}_h^{n}\cdot\mathbf{n}\hat{c}_h^{n},\mu_h-v_h\rangle_{\partial\mathcal{T}_h}.
\end{aligned}
\]
Now, we can arrive at the following hybrid mixed DG finite element method.
\begin{alg}[DGFE Algorithm]\label{alg4}
 For given $\mathbf{u}_h^{n}$, find $(\bm\sigma_h^{n},c_h^{n},\lambda_h^{c,n})\in \Theta_h \times \Lambda_h \times \Sigma_h$ such that
\begin{equation}\label{hmdg}
\begin{aligned}
B_c(\mathbf{u}^{n}_h;(\bm\sigma_h^{n},c_h^{n},&\lambda_h^{c,n}),(\bm\tau_h,v_h,\mu_h))=-(q^{n}c^{*,n}_{h}+\phi c_h^{n-1}/ \Delta t,v_h)_{\mathcal{T}_h},
\\& \forall (\bm\tau_h,v_h,\mu_h) \in \Theta_h \times \Lambda_h \times \Sigma_h.
\end{aligned}
\end{equation}
\end{alg}
\subsection{The hybrid mixed discontinuous Galerkin finite element method}

Now we give the combined hybrid mixed discontinuous Galerkin finite element method for incompressible miscible displacement problem.
\begin{alg}[HMDGFE Algorithm]\label{alg5}
For the given initial approximate values $c^{0}_{h}$, find $(\mathbf{u}_h^{n}, p_h^{n}, \lambda_h^{n}) \in \Theta_h \times \Lambda_h \times \Sigma_h$ and $(\bm\sigma_h^{n},c_h^{n},\lambda_h^{c,n})\in \Theta_h \times \Lambda_h \times \Sigma_h$, such that, for $ \forall (\bm\tau_h,v_h,\mu_h) \in \Theta_h \times \Lambda_h \times \Sigma_h$
\begin{equation}\label{101}
\begin{aligned}
&({\rm a})\quad B_{\mathbf{u}}(c_h^{n-1};(\mathbf{u}_h^{n},p_h^{n},\lambda_h^{n}),(\bm\tau_h,v_h,\mu_h))=-(q^{n},v_h)_{\mathcal{T}_h},\\
&({\rm b})\quad B_c(\mathbf{u}^{n}_h;(\bm\sigma_h^{n},c_h^{n},\lambda_h^{c,n}),(\bm\tau_h,v_h,\mu_h))=-(q^{n}c^{*,n}+\frac{\phi}{\Delta t}c_h^{n-1},v_h)_{\mathcal{T}_h}.
\end{aligned}
\end{equation}
\end{alg}
\begin{theorem}[Consistency]
Algorithm \ref{alg5} is consistent. That is, let $\lambda = p$, $\lambda^{c} = c$ and $\bm\sigma=-D(\mathbf{u})\nabla c$ and $(\mathbf{u},p,\bm\sigma,c)$ be the solution of the problem \eqref{eq1}, then the corresponding variational equation \eqref{101} holds if $(\mathbf{u}_h^{n},p_h^{n},\lambda_{h}^{n},\bm\sigma_h^{n},c_h^{n},\lambda_h^{c,n})$ are replaced by
$(\mathbf{u}^{n},p^{n},\lambda^{n},\bm\sigma^{n},c^{n},\lambda^{c,n})$ for $n=0,1,2,\ldots,M$.
\end{theorem}
\begin{proof}
 Substituting the solution $(\mathbf{u}^{n},p^{n})$ of problem \eqref{eq1} into \eqref{101}(a) with $v_h=\mu_h=0$ and $\bm\tau_h=\mu_h=0$ respectively, we can get
 \[
B_\mathbf{u}(c^{n-1};(\mathbf{u}^{n},p^{n},p^{n});(\bm\tau_h,0,0))=(a(c^{n-1})\mathbf{u}^{n},\bm\tau_h)_{\mathcal{T}_h}+(\nabla p^{n},\bm\tau_h)_{\mathcal{T}_h}=0
\]
and
\[
\begin{aligned}
B_\mathbf{u}(c^{n-1};(\mathbf{u}^{n},p^{n},p^{n});(0,v_h,0))=(\mathbf{u}^{n},\nabla v_h)_{\mathcal{T}_h}-\langle\mathbf{u}^{n}\cdot\mathbf{n},v_h\rangle_{\partial\mathcal{T}_h}
=-(\nabla\cdot\mathbf{u}^{n},v_h)_{\mathcal{T}_h}
=-(q^{n},v_h)_{\mathcal{T}_h}.
\end{aligned}
\]

Next we test with $\bm\tau_h=v_h= 0$ and get the equation
\[
\begin{aligned}
B_\mathbf{u}(c^{n-1};(\mathbf{u}^{n},p^{n},p^{n});(0,0,\mu_h))=\langle\mathbf{u}^{n}\cdot\mathbf{n},\mu_h\rangle_{\partial\mathcal{T}_h}=0,
\end{aligned}
\]
thus the normal flux $\mathbf{u}^{n}\cdot\mathbf{n}$ is continuous across element interfaces.

Now, we consider HMDG algorithm \eqref{101}(b) and let $(\bm\sigma^{n},c^{n})$ denote the solution of \eqref{eq1}. Now testing with $(0,v_h,0)$ and $(0,0,\mu_h)$ respectively, we have
\[
\begin{aligned}
&B_c(\mathbf{u}^{n};(\bm\sigma^{n},c^{n},c^{n});(0,v_h,0))\\
=&-(\frac{\phi}{\Delta t}c^{n},v_h)_{\mathcal{T}_h}+(\bm\sigma^{n}+\mathbf{u}^{n}c^{n},\nabla v_h)_{\mathcal{T}_h}-\langle\bm\sigma^{n}\cdot\mathbf{n}+\mathbf{u}^{n}\cdot\mathbf{n}c^{n},v_h\rangle_{\partial\mathcal{T}_h}\\
=&-(\frac{\phi}{\Delta t}c^{n},v_h)_{\mathcal{T}_h}-(\nabla\cdot(-D(\mathbf{u}^{n})\nabla c^{n}+\mathbf{u}^{n}c^{n}),v_h)_{\mathcal{T}_h}\\
=&-(q^{n}c^{*,n}+\frac{\phi}{\Delta t}c^n,v_h)_{\mathcal{T}_h},
\end{aligned}
\]
and
\[
B_c(\mathbf{u}^{n};(\bm\sigma^{n},c^{n},c^{n});(0,0,\mu_h))=\langle\bm\sigma^{n}\cdot\mathbf{n}+\mathbf{u}^{n}\cdot\mathbf{n}c^{n},\mu_h\rangle_{\partial\mathcal{T}_h}=0.
\]
Thus we have proved the consistency of Algorithm \ref{alg5}.
\end{proof}
\begin{theorem}[Mass conservation]\label{conservation}

Algorithm \ref{alg5} is locally and globally conservative.
\end{theorem}
\begin{proof}

First, we will show the local conservation of HMFE scheme \eqref{101}(a) with the test functions $(0,1,0)$ on $K$. This yields
\[
\begin{aligned}
B_\mathbf{u}(c_h^{n-1};(\mathbf{u}_h^{n},p_h^{n},\lambda_h^{n});(0,1,0))=-\langle\mathbf{u}_h^{n}\cdot\mathbf{n},1\rangle_{\partial K}=-(q^{n},1)_{K}.\\
\end{aligned}
\]
Thus the total flux over an element boundary equals the sum of internal sources, and hence HMFE scheme \eqref{101}(a) is locally conservative. We can obtain the continuity of the normal fluxes $\mathbf{u}_h\cdot\mathbf{n}$ across element interfaces with $(0,0,1)$ for some $e\in\partial\mathcal{T}_h$, so HMFE scheme \eqref{101}(a) is also globally conservative.

Next,  as in the proof of HMFE scheme \eqref{101}(a), we show the local conservation of HMDGFE scheme \eqref{101}(b). With $(0,1,0)$ on $K$, we can get
\[
\begin{aligned}
&B_c(\mathbf{u}_h^{n};(\bm\sigma_h^{n},c_h^{n},\lambda_h^{c,n});(0,1,0))\\
=&-(\frac{\phi}{\Delta t}c_h^{n},1)_{K}-\langle\bm\sigma_h^{n}\cdot\mathbf{n},1\rangle_{\partial K}-\langle\mathbf{u}_h^{n}\cdot\mathbf{n}c_h^{n},1\rangle_{\partial K^{out}}
-\langle\mathbf{u}_h^{n}\cdot\mathbf{n}\lambda_h^{c,n},1\rangle_{\partial K^{in}}\\
=&-(q^{n}c^{*,n}+\frac{\phi}{\Delta t}c_h^{n-1},1)_{K}.
\end{aligned}
\]
So the total flux over the element boundaries equals the sum of internal sources and the flux over the boundary of the domain. Now let $e\in\partial\mathcal{T}_h$ and $e=\partial K_1^{in}\bigcap\partial K_2^{out}$. Testing with $(0,0,1)$ on $e$, we obtain as follows
\begin{equation}\nonumber
\begin{aligned}
B_c(\mathbf{u}_h^{n};(\bm\sigma_h^{n},c_h^{n},\lambda_h^{c,{n}});(0,0,1))
=\langle\bm\sigma_h^{n}\cdot\mathbf{n},1\rangle_{e}+\langle\mathbf{u}_h^{n}\cdot\mathbf{n}c_h,1\rangle_{\partial K_2^{out}}+\langle\mathbf{u}_h^{n}\cdot\mathbf{n}\lambda_h^{cn+1},1\rangle_{\partial K_1^{in}}=0,
\end{aligned}
\end{equation}
so we know that the total outflow over a facet on one element balances the inflow over the same facet on the neighbouring element.
\end{proof}

\section{Some preliminaries and convergence theorem}

In order to prove the convergence of our proposed algorithm, we will use the following important result in \cite{hmixed1}.
\begin{lemma}\label{test_func}
There is a unique solution $\tilde{\tau}_c\in \Theta_h$ defined elementwise by the variational problem
\[
\begin{aligned}
(\tilde{\tau}_c,\omega)_K&=(\nabla c_h,\omega)_K,\quad\quad\forall \omega\in[P_{k-1}(K)]^d,\\
\langle\tilde{\tau}_c\cdot\mathbf{n},\mu\rangle_{\partial K}&=\langle\frac{1}{h}(\lambda^{c}_h-c_h),\mu\rangle_{\partial K},\quad\quad \mu\in P_k(\partial K),
\end{aligned}
\]
where  $\lambda^{c}_h\in\Sigma_h$.

Moreover, the following result
\[
\begin{aligned}
\|\tilde{\tau}_c\|_{\mathcal{T}_h}\leq C_c(\|\nabla c_h\|_{\mathcal{T}_h}^2+\frac{1}{h}|\lambda^{c}_h-c_h|_{\partial\mathcal{T}_h}^2)^{\frac{1}{2}}
\end{aligned}
\]
holds, where  $C_c$  is  a constant independent of the parameters $h$ and $\Delta t$.
\end{lemma}

The uniformly positive definiteness and Lipschitz continuousness of the diffusion coefficient $D(\mathbf{u})$ as in \cite{hgmixed3} will be used.
\begin{lemma}\label{DU}
If the molecular diffusion and dispersion coefficients $d_m$, $d_l$ and $d_t$ are nonnegative, then we have
\begin{equation}
\begin{aligned}\nonumber
(a)\quad&D(\mathbf{u})^{-1}\mathbf{v}\cdot\mathbf{v}\leq(d_m+\min(d_l,d_t)| \mathbf{u}|)^{-1}|\mathbf{v}|^2\leq d_m^{-1}|\mathbf{v}|^2,\\
(b)\quad&D(\mathbf{u})^{-1}\mathbf{v}\cdot\mathbf{v}\geq(d_m+max(d_l,d_t)| \mathbf{u}|)^{-1}|\mathbf{v}|^2,\\
(c)\quad&|D(\mathbf{u})^{-1}-D(\mathbf{v})^{-1}|\leq d_m^{-2}(7d_t+6d_l)d^{\frac{3}{2}}|\mathbf{u}-\mathbf{v}|.
\end{aligned}
\end{equation}
\end{lemma}

Using the similar techniques as in \cite{hmixed1,hmixed2,hgmixed1,hgmixed2}, we can easily get the following stability and boundedness of the bilinear form $B_\mathbf{u}$ under the pair of the norms
\[
\begin{aligned}
&||(\mathbf{u},p,\lambda)||_\mathbf{u}:=(||\mathbf{u}||_{\mathcal{T}_h}^2+||\nabla p||_{\mathcal{T}_h}^2+\frac{1}{h}|\lambda-p|_{\partial\mathcal{T}_h}^2)^\frac{1}{2},\\
&||(\mathbf{u},p,\lambda)||_\mathbf{u,*}:=(||(\mathbf{u},p,\lambda)||^{2}_\mathbf{u}+h|\mathbf{u}\cdot\mathbf{n}|^{2}_{\p \mathcal{T}_{h}})^\frac{1}{2}.
\end{aligned}
\]
\begin{lemma}[Stability and boundedness of $B_{\mathbf{u}}$]\label{inf_sup_u}
There exist two positive constants $C^{\mathbf{u},\ast}$ and $C_{\mathbf{u},\ast}$
that are independent of the mesh size $h$ such that
\begin{equation}\nonumber
\begin{aligned}
&(\up{a})\quad| B_{\mathbf{u}}(c;(\mathbf{u}, p, \lambda),(\bm\tau_h, v_h, \mu_h))|\leq C^{\mathbf{u},\ast}\|(\mathbf{u}, p, \lambda)\|_{\mathbf{u,*}}\|(\bm\tau_h, v_h, \mu_h)\|_\mathbf{u},\\
&(\up{b})\quad\sup\limits_{(\bm\tau_h, v_h, \mu_h)}\frac{B_{\mathbf{u}}(c_h;(\mathbf{u}_h, p_h, \lambda_h),(\bm\tau_h, v_h, \mu_h))}{\|(\bm\tau_h, v_h, \mu_h)\|_\mathbf{u}}\geq C_{\mathbf{u},\ast}\|(\mathbf{u}_h, p_h, \lambda_h)\|_\mathbf{u},
\end{aligned}
\end{equation}
holds for all $( \mathbf{u}_h, p_h, \lambda_h) \in \Theta_h \times \Lambda_h \times \Sigma_h$ and $(\bm\tau_h, v_h, \mu_h) \in \Theta_h \times \Lambda_h \times \Sigma_h$.
\end{lemma}
To show the stability and boundedness of the bilinear form $B_c$, we definite the norms as follows
\begin{equation}
\begin{aligned}
||(\bm\sigma,c,\lambda^c)||_{D}&:=(||\bm\sigma||_{\mathcal{T}_h}^2+\f{1}{\Delta t}|| c||_{\mathcal{T}_h}^2+||\nabla c||_{\mathcal{T}_h}^2+\frac{1}{h}|\lambda^c-c|_{\partial\mathcal{T}_h}^2)^\frac{1}{2},\\
||(c,\lambda^c)||_B&:=(||\mathbf{u}\cdot\nabla c||_{\mathcal{T}_h}^2+||\mathbf{u}\cdot\mathbf{n}|^{1/2}(\lambda^c-c)|_{\partial\mathcal{T}_h}^2)^\frac{1}{2},\\
||(\bm\sigma,c,\lambda^c)||_{B,*}&:=(h|\bm\sigma\cdot\mathbf{n}|_{\p\mathcal{T}_h}^2+\|\mathbf{u}\cdot\nabla c\|_{\partial\mathcal{T}_h}^2+|\mathbf{u}\cdot\mathbf{n}\lambda^{c}|^{2}_{\p\mathcal{T}_{h}})^\frac{1}{2},\\
||(\bm\sigma,c,\lambda^c)||_c&:=(||(\bm\sigma,c,\lambda^c)||_D^2+||(c,\lambda^c)||_B^2)^{\frac{1}{2}},
\\
||(\bm\sigma,c,\lambda^c)||_{c,*}&:=(||(\bm\sigma,c,\lambda^c)||_D^2+||(c,\lambda^c)||_{B,*}^2)^{\frac{1}{2}}.
\end{aligned}
\end{equation}
We have the following the stability and boundedness result on  the bilinear form $B_c$.


\begin{lemma}[Stability and boundedness of $B_{c}$]\label{inf_sup_c}
There exist two positive constants $C^{c,\ast}$ and $C_{c,\ast}$
that are independent of the mesh size $h$ such that, for some given $\Delta t_{0}>0$, when $\Delta t\leq \Delta t_{0}$,
\begin{equation}\nonumber
\begin{aligned}
&(\up{a})\quad| B_c(\mathbf{u};(\bm\sigma,c,\lambda^{c}),(\bm\tau_h,v_h,\mu_h))|\leq C^{c,\ast}\|(\bm\sigma,c,\lambda^{c})\|_{c,*}\|(\bm\tau_h,v_h,\mu_h)\|_c,\\
&(\up{b})\quad\sup\limits_{(\bm\tau_h,v_h,\mu_h)}\frac{B_c(\mathbf{u}_h;(\bm\sigma_h,c_h,\lambda_h^{c}),(\bm\tau_h,v_h,\mu_h))}{\|(\bm\tau_h, v_h, \mu_h)\|_c}\geq C_{c,\ast}\|(\bm\sigma_h,c_h,\lambda_h^{c})\|_c,\\
\end{aligned}
\end{equation}
holds for all $( \bm\sigma_h,c_h,\lambda_h^{c}) \in \Theta_h \times \Lambda_h \times \Sigma_h$ and $(\bm\tau_h, v_h, \mu_h) \in \Theta_h \times \Lambda_h \times \Sigma_h$.
\end{lemma}
\begin{proof}
Firstly, we choose $(\bm\tau_{h},v_{h},\mu_{h})=(\gamma\tilde{\bm\tau}_{c},0,0)$ in the bilinear form $B_{c}$ and use Lemma \ref{test_func} to get
\begin{equation}\label{eq3-1}
\begin{split}
&B_{c}(\mathbf{u}_{h};(\bm\sigma_h,c_h,\lambda_h^{c}),(\gamma\tilde{\bm\tau}_{c},0,0)
\\=&(\frac{1}{D(\mathbf{u}_h)}\bm\sigma_h,\gamma\tilde{\bm\tau}_{c})_{\mathcal{T}_h}+(\nabla c_{h},\gamma\tilde{\bm\tau}_{c})_{\mathcal{T}_h}+\langle\lambda_h^{c}-c_h,\gamma\tilde{\bm\tau}_{c}\cdot\mathbf{n}\rangle_{\partial\mathcal{T}_h}
\\=&\gamma(\frac{1}{D(\mathbf{u}_h)}\bm\sigma_h,\tilde{\bm\tau}_{c})_{\mathcal{T}_h}+\gamma\|\nabla c_h\|^{2}_{\mathcal{T}_h}+\f{\gamma}{h}|\lambda_h^{c}-c_h|^{2}_{\partial\mathcal{T}_h}
\\\geq&-\f{1}{2}(\frac{1}{D(\mathbf{u}_h)}\bm\sigma_h,\bm\sigma_h)_{\mathcal{T}_h}-\f{\gamma^{2}}{2d_m}\|\tilde{\bm\tau}_{c}\|^{2}_{\mathcal{T}_h}+\gamma(\|\nabla c_h\|^{2}_{\mathcal{T}_h}+\f{1}{h}|\lambda_h^{c}-c_h|^{2}_{\partial\mathcal{T}_h})
\\\geq&-\f{1}{2}(\frac{1}{D(\mathbf{u}_h)}\bm\sigma_h,\bm\sigma_h)_{\mathcal{T}_h}+(\gamma-\f{C^{2}_{c}\gamma^{2}}{2d_m})(\|\nabla c_h\|^{2}_{\mathcal{T}_h}+\f{1}{h}|\lambda_h^{c}-c_h|^{2}_{\partial\mathcal{T}_h}).
\end{split}
\end{equation}

And then, taking $(\bm\tau_{h},v_{h},\mu_{h})=(\bm\sigma_{h},-c_{h},-\lambda^{c}_{h})$  in the bilinear form $B_{c}$, we have
\begin{equation}\label{eq3-2}
\begin{split}
&B_{c}(\mathbf{u}_{h};(\bm\sigma_h,c_h,\lambda_h^{c}),(\bm\sigma_{h},-c_{h},-\lambda^{c}_{h})
\\=&\f{1}{\Delta t}(\phi c_h,c_h)_{\mathcal{T}_h}+(\frac{1}{D(\mathbf{u}_h)}\bm\sigma_h,\bm\sigma_h)_{\mathcal{T}_h}
-(\mathbf{u}_hc_h,\nabla c_h)_{\mathcal{T}_h}+\langle\mathbf{u}_h\cdot\mathbf{n}\hat{c}_h,-\lambda_h^{c}+c_h\rangle_{\partial\mathcal{T}_h}
\\=&\f{1}{\Delta t}(\phi c_h,c_h)_{\mathcal{T}_h}+\f{1}{2}(\nabla\cdot\mathbf{u}_{h} c_h,c_h)_{\mathcal{T}_h}+(\frac{1}{D(\mathbf{u}_h)}\bm\sigma_h,\bm\sigma_h)_{\mathcal{T}_h}
\\&-\f{1}{2}\langle\mathbf{u}_h\cdot\mathbf{n}{c}_h,c_h\rangle_{\partial\mathcal{T}_h}+\langle\mathbf{u}_h\cdot\mathbf{n}\hat{c}_h,-\lambda_h^{c}+c_h\rangle_{\partial\mathcal{T}_h}
\\=&\f{1}{\Delta t}(\phi c_h,c_h)_{\mathcal{T}_h}+\f{1}{2}(q c_h,c_h)_{\mathcal{T}_h}+(\frac{1}{D(\mathbf{u}_h)}\bm\sigma_h,\bm\sigma_h)_{\mathcal{T}_h}
\\&-\f{1}{2}\langle\mathbf{u}_h\cdot\mathbf{n}{c}_h,c_h\rangle_{\partial\mathcal{T}_h}+\langle\mathbf{u}_h\cdot\mathbf{n}\hat{c}_h,-\lambda_h^{c}+c_h\rangle_{\partial\mathcal{T}_h}.
\end{split}
\end{equation}
Note that
\begin{equation}\label{eq3-3}
\begin{aligned}
&(\up{a})\quad-\frac{1}{2}\langle\mathbf{u}_h\cdot\mathbf{n}c_h,c_h\rangle_{\partial\mathcal{T}_h}=\frac{1}{2}\langle|\mathbf{u}_h\cdot\mathbf{n}|c_h,c_{h}\rangle_{\partial\mathcal{T}_h^{in}}-\frac{1}{2}\langle|\mathbf{u}_h\cdot\mathbf{n}|c_h,c_{h}\rangle_{\partial\mathcal{T}_h^{out}},\\
&(\up{b})\quad-\langle\mathbf{u}_h\cdot\mathbf{n}\hat{c}_h,\lambda_h^{c}\rangle_{\partial\mathcal{T}_h}=\langle|\mathbf{u}_h\cdot\mathbf{n}|\lambda_h^{c},\lambda_h^{c}\rangle_{\partial\mathcal{T}_h^{in}}-\langle|\mathbf{u}_h\cdot\mathbf{n}| \lambda_h^{c},c_h\rangle_{\partial\mathcal{T}_h^{out}},\\
&(\up{c})\quad\langle\mathbf{u}_h\cdot\mathbf{n}\hat{c}_h,c_h\rangle_{\partial\mathcal{T}_h}=\langle|\mathbf{u}_h\cdot\mathbf{n}|c_h,c_h\rangle_{\partial\mathcal{T}_h^{out}}-\langle |\mathbf{u}_h\cdot\mathbf{n}|\lambda_h^{c},c_h\rangle_{\partial\mathcal{T}_h^{in}}.
\end{aligned}
\end{equation}
Now let $K_1$ and $K_2$ denote two elements sharing the facet  $e=\partial K_1^{out}\bigcap K_2^{in}$. Since $\lambda_h^{c}$ is a single value function on $e$, we have $\lambda_h^{c}|_{\partial K_1^{out}}=\lambda_h^{c}|_{\partial K_2^{in}}$, which means that we can shift the terms only involving the Lagrange multiplier between neighbouring elements. Hence,  summing \eqref{eq3-3} up, we can  rewrite the last two terms of \eqref{eq3-2} as follows
\[
\frac{1}{2}\langle |\mathbf{u}_h\cdot\mathbf{n}|(\lambda_h^{c}-c_{h}),\lambda_h^{c}-c_h\rangle_{\partial\mathcal{T}_h}.
\]
And then, we get
\begin{equation}\label{eq3-4}
\begin{aligned}
&B_{c}(\mathbf{u}_{h};(\bm\sigma_h,c_h,\lambda_h^{c}),(\bm\sigma_{h},-c_{h},-\lambda^{c}_{h})
\\=&((\f{\phi}{\Delta t}+\f{q}{2}) c_h,c_h)_{\mathcal{T}_h}+(\frac{1}{D(\mathbf{u}_h)}\bm\sigma_h,\bm\sigma_h)_{\mathcal{T}_h}+\frac{1}{2}\langle |\mathbf{u}_h\cdot\mathbf{n}|(\lambda_h^{c}-c_{h}),\lambda_h^{c}-c_h\rangle_{\partial\mathcal{T}_h}.
\end{aligned}
\end{equation}

In \eqref{eq3-4}, we choose some time step $\Delta t_{0}>0$ such that $\f{\phi}{\Delta t}+\f{q}{2}\geq C_{*}>0$ when $\Delta t\leq \Delta t_{0}$,  and we also take $\gamma =d_{m}/C^{2}_{c}$ ($\gamma-\f{C^{2}_{c}\gamma^{2}}{2d_m}=\f{d_{m}}{2C^{2}_{c}}>0$) in \eqref{eq3-1}. Then combing \eqref{eq3-1} and \eqref{eq3-4}, we obtain the stability of the bilinear form $B_{c}$. Using Cauchy inequality we can easily get the boundedness of the bilinear form $B_{c}$. Here we are not going to give any details.

\end{proof}

By the stability and boundedness of the bilinear forms $B_{D}$ and $B_{c}$, with Lax-Milgram theorem, the following existence theorem can be easily obtained.
\begin{theorem}[Existence and Uniqueness]
For given initial approximate value $c^{0}_{h}$, there exists a parameter $\Delta t_{0}>0$, such that, when $\Delta  t\leq \Delta t_{0}$, HMDGFE Algorithm exists a unique solution.
\end{theorem}

Next, we will give some important projection operators and approximate properties, which is used to show the convergence theorem of our proposed method.

Introduce the local $L^2$-projection operators
$\Pi_h$ and $\Pi_e$ as follows:
\begin{equation}\nonumber
\begin{aligned}
&(p-\Pi_hp,v_h)_K=0,\quad \forall v_h\in P_k(K),\\
&\langle\lambda-\Pi_e\lambda,\mu_h\rangle_e=0,\quad \forall \mu_h\in P_k(e),
\end{aligned}
\end{equation}
 where $K\in\mathcal{T}_h$, $e\in\partial\mathcal{T}_h$, $p\in L^2(K)$ and $\lambda\in L^2(e)$. These projection operators satisfy the following error estimates \cite{ref21}.
\begin{lemma}\label{L2_err}
For the local $L^2$-projection operators $\Pi_h$ and $\Pi_e$, we have the estimate
\begin{equation}\nonumber
\begin{aligned}
&||p-\Pi_hp||_K\leq Ch^s||p||_{s,K},\quad 0\leq s\leq k+1,\\
&||\nabla(p-\Pi_hp)||_K\leq Ch^s||p||_{s+1,K},\quad 0\leq s\leq k,\\
&|p-\Pi_hp|_e+|p-\Pi_ep|_e\leq Ch^{s+\frac{1}{2}}||p||_{s+1,K},\quad 0\leq s\leq k,
\end{aligned}
\end{equation}
 where $C$ is a constant independent of $h$.
\end{lemma}
\begin{lemma}\label{pro_err_c}
Suppose that Assumption $1.2$ holds. Then, for any element $K$, we can reach
\begin{equation}\nonumber
\begin{aligned}
&||\nabla(c-\Pi_hc)||_K\leq Ch^s||c||_{s+1,K},\quad 0\leq s\leq k+1,\\
&|c-\Pi_ec|_{\partial K}\leq Ch^{s}|c|_{s,\partial K},\quad 0\leq s\leq k+1,
\end{aligned}
\end{equation}
where $C$ is a constant independent of  $h$.
\end{lemma}
Similarly, the interpolation operators for functions on $\mathcal{T}_h$ and $\partial\mathcal{T}_h$ are defined element-wise and are denoted by the same symbols.
For $\mathbf{u}\in H(\up{div},K)$ we utilize the Raviart-Thomas interpolantion \cite{hmixed1} defined by
\begin{equation}\nonumber
\begin{aligned}
&(\mathbf{u}-\Pi^{RT}\mathbf{u},\omega_h)_K=0,\quad\forall\omega_h\in[P_{k-1}(K)]^d,\\
&\langle(\mathbf{u}-\Pi^{RT}\mathbf{u})\cdot \mbf{n}_e,\mu_h\rangle_e=0,\quad\forall\mu_h\in P_{k}(e),\quad e\in\partial K.
\end{aligned}
\end{equation}
We can reach the error estimate as follows:
\begin{lemma}\label{rt}
For the projection operator $\Pi^{RT}$ defined as above, we have the estimate
\begin{equation}\nonumber
\begin{aligned}
&||\mathbf{u}-\Pi^{RT}\mathbf{u}||_K+h^{\frac{1}{2}}|u-\Pi^{RT}\mathbf{u}|_{\partial K}\leq Ch^s||\mathbf{u}||_{s,K},\quad \frac{1}{2}\leq s\leq k+1,\\
&||\nabla\cdot(\mathbf{u}-\Pi^{RT}\mathbf{u})||_K\leq Ch^s||\nabla\cdot\mathbf{u}||_{s,K},\quad 1\leq s\leq k+1,
\end{aligned}
\end{equation}
where $C$ is a constant independent of $h$.
\end{lemma}

The following trace inequalities will be also used to prove the convergence theorem.
\begin{lemma}
For $\forall v\in H^1(K)$, the trace inequalities are shown below
\begin{equation}\label{1111}
\begin{aligned}
&\| v\|_{0,e}^2\leq C(h_e^{-1}\| v\|_{0,K}^2+h_e\| v\|_{1,K}^2),\\
&\| \nabla v\cdot \mbf{n}_e\|_{0,e}^2\leq C(h_e^{-1}\| \nabla v\|_{0,K}^2+h_e\|\nabla^2v\|_{0,K}^2).
\end{aligned}
\end{equation}
\end{lemma}

For HMDGFE Algorithm, we have the following main convergence theorem.
\begin{theorem}[Convergence theorem]\label{err_c}
Suppose that  Assumptions $1.1$ and $1.2$ hold. And let $(\mathbf{u}_h, p_{h}, \lambda_{h},\bm\sigma_h,c_h,\lambda^{c}_{h})$ be the solution of HMDGFE Algorithm with initial values $c^0=c_h^{0}=\Pi_hc^{0}$. Then, for some given $\Delta t_{0}>0$, when $\Delta t\leq \Delta t_{0}$, we have the following error estimate,  for $m>0$
\begin{equation}\label{eq24}
\begin{aligned}
(\up{a})\quad&\|\mathbf{u}^{m}-\mathbf{u}_h^{m}\|_{\mathcal{T}_h}^2+\|p^{m}-p_h^{m}\|_{\mathcal{T}_h}^2\leq C(h^{2s}+\Delta t^2),\\
(\up{b})\quad&\|c^{m}-c_h^{m}\|_{\mathcal{T}_h}^2+\Delta t\sum\limits_{n=0}^m\|\bm\sigma^{n}-\bm\sigma_h^{n}\|_{\mathcal{T}_h}^2\leq C(h^{2s}+\Delta t^2),\\
(\up{c})\quad&\|\nabla(\Pi_hp^{m}-p_h^{m})\|_{\mathcal{T}_h}^2+\frac{1}{h}|\lambda_h^{m}-p_h^{m}|_{\partial\mathcal{T}_h}^2\leq C(h^{2s}+\Delta t^2),\\
(\up{d})\quad&\Delta t\sum\limits_{n=0}^mh|\lambda_h^{c,n}-c^{n}|_{\partial\mathcal{T}_h}^2\leq C(h^{2s}+\Delta t^2),
\end{aligned}
\end{equation}
where $C$ is a constant independent of $h$ and $\Delta t$, and $0 \leq s \leq k + 1$ when $d = 2$, $1 \leq s \leq k + 1$ when $d = 3$.
\end{theorem}

\section{Proof of convergence theorem}

For HMFE Algorithm, we have the following  error estimate:
\begin{lemma}\label{err_u}
Suppose that the coefficients of system satisfy Assumption \ref{eq3}. Then, for any $n>0$, the following inequality holds:
\begin{equation}\nonumber
\begin{aligned}
&\|(\mathbf{u}^{n}-\mathbf{u}_h^{n},\Pi_{h}p^{n}-p_h^{n},\Pi_{e}p^{n}-\lambda_h^{n})\|_{\mathbf{u}}
\leq C(\|\Pi^{RT}\mathbf{u}^{n}-\mathbf{u}^{n}\|_{\mathcal{T}_h}+\|c^{n-1}-c_h^{n-1}\|_{\mathcal{T}_h}+\tr t)
\end{aligned}
\end{equation}
where $C$ is a constant independent of $h$ and $\Delta t$.
\end{lemma}
\begin{re}
Using the similar technique as in  \cite{hgmixed1}, we can easily get the proof of Lemma \ref{err_u}.
\end{re}

\begin{lemma}\label{err_gradc}
 Under Assumption \ref{eq3}, for any $m>0$, the following inequality holds:
\begin{equation}\nonumber
\begin{aligned}
\|\Pi_hc^m-c_h^m\|_{\mathcal{T}_h}^2+\Delta t\sum_{n=1}^m\|\nabla(\Pi_hc^n-c_h^n)\|_{\mathcal{T}_h}^2\leq C\Delta t\sum_{n=1}^m\|\Pi^{RT}\bm\sigma^n-\bm\sigma_h^n\|_{\mathcal{T}_h}^2
\end{aligned}
\end{equation}
where $C$ is a constant independent of $h$ and $\Delta t$.
\end{lemma}
\begin{proof}
First, we denote the following notation for convenience
\begin{equation}
\begin{aligned}
&\zeta_c=\Pi_hc-c_h,\quad\zeta_{\sigma}=\Pi^{RT}\bm\sigma-\bm\sigma_h,\quad\zeta_\lambda=\Pi_ec-\lambda^c_h.
\end{aligned}
\end{equation}
We can derive the following formulas from the consistency of DGFE Algorithm
\begin{equation}
\begin{aligned}
&B_c(\mathbf{u}_h^{n};(\zeta_\sigma^n,\zeta_c^n,\zeta_\lambda^n);(\bm\tau_h,v_h,\mu_h))\\
=&-\frac{1}{\Delta t}(\phi\zeta_c^n,v_h)_{\mathcal{T}_h}+(\zeta_\sigma^n,\nabla v_h)_{\mathcal{T}_h}+(\mathbf{u}_h^{n}\zeta_c^n,\nabla v_h)_{\mathcal{T}_h}
+(D(\mathbf{u}_h^{n})^{-1}\zeta_\sigma^n,\bm\tau_h)_{\mathcal{T}_h}\\
&+(\nabla\zeta_c^n,\bm\tau_h)_{\mathcal{T}_h}+\langle\zeta_\lambda^n,\bm\tau_h\cdot\mathbf{n}\rangle_{\partial\mathcal{T}_h}
-\langle\zeta_c^n,\bm\tau_h\cdot\mathbf{n}\rangle_{\partial\mathcal{T}_h}+\langle\zeta_\sigma^n\cdot\mathbf{n},\mu_h-v_h\rangle_{\partial\mathcal{T}_h}\\
&+\langle\mathbf{u}_h^{n}\cdot\mathbf{n}\zeta_\lambda^n,\mu_h-v_h\rangle_{\partial\mathcal{T}_h^{in}}+\langle\mathbf{u}_h^{n}\cdot\mathbf{n}\zeta_c^n,\mu_h-v_h\rangle_{\partial\mathcal{T}_h^{out}}\\
=&-(q\zeta_c^{n,\ast},v_h)_{\mathcal{T}_h}-\frac{1}{\Delta t}(\phi\zeta_c^{n-1},v_h)_{\mathcal{T}_h}.
\end{aligned}
\end{equation}
Setting $v_h=-\zeta_c^n$ and $\mu_h=-\zeta_\lambda^n$ and choosing $\bm\tau_h=\nabla\zeta_c^n$ for every element $K\in\mathcal{T}_h$ and $\bm\tau_h\cdot\mathbf{n}=\frac{1}{h}(\zeta_\lambda^n-\zeta_c^n)$ for $e\in\partial\mathcal{T}_h$, so we can obtain
\begin{align*}
&B_c(\mathbf{u}_h^{n};(\zeta_\sigma^n,\zeta_c^n,\zeta_\lambda^n);(\bm\tau_h,v_h,\mu_h))\\
=&\frac{1}{\Delta t}(\phi\zeta_c^n,\zeta_c^n)_{\mathcal{T}_h}-(\zeta_\sigma^n,\nabla \zeta_c^n)_{\mathcal{T}_h}-(\mathbf{u}_h^{n}\zeta_c^n,\nabla \zeta_c^n)_{\mathcal{T}_h}
+(D(\mathbf{u}_h^{n})^{-1}\zeta_\sigma^n,\nabla \zeta_c^n)_{\mathcal{T}_h}
+(\nabla\zeta_c^n,\nabla\zeta_c^n)_{\mathcal{T}_h}\\
&+\langle\zeta_\lambda^n,\frac{1}{h}(\zeta_\lambda^n-\zeta_c^n)\rangle_{\partial\mathcal{T}_h}
-\langle\zeta_c^n,\frac{1}{h}(\zeta_\lambda^n-\zeta_c^n)\rangle_{\partial\mathcal{T}_h}+\langle\zeta_\sigma^n\cdot\mathbf{n},\zeta_c^n-\zeta_\lambda^n\rangle_{\partial\mathcal{T}_h}\\
&+\langle\mathbf{u}_h^{n}\cdot\mathbf{n}\zeta_\lambda^n,\zeta_c^n-\zeta_\lambda^n\rangle_{\partial\mathcal{T}_h^{in}}+\langle\mathbf{u}_h^{n}\cdot\mathbf{n}\zeta_c^n,\zeta_c^n-\zeta_\lambda^n\rangle_{\partial\mathcal{T}_h^{out}}\\
=&\frac{1}{\Delta t}(\phi\zeta_c^n,\zeta_c^n)_{\mathcal{T}_h}-(\zeta_\sigma^n,\nabla \zeta_c^n)_{\mathcal{T}_h}+\frac{1}{2}(\nabla\cdot\mathbf{u}_h^{n}\zeta_c^n,\zeta_c^n)_{\mathcal{T}_h}
+(D(\mathbf{u}_h^{n})^{-1}\zeta_\sigma^n,\nabla \zeta_c^n)_{\mathcal{T}_h}\\
&+(\nabla\zeta_c^n,\nabla\zeta_c^n)_{\mathcal{T}_h}-\frac{1}{2}\langle\mathbf{u}_h^{n}\cdot\mathbf{n}\zeta_c^n,\zeta_c^n\rangle_{\partial\mathcal{T}_h}
+\langle\zeta_\lambda^n-\zeta_c^n,\frac{1}{h}(\zeta_\lambda^n-\zeta_c^n)\rangle_{\partial\mathcal{T}_h}\\
&+\langle\zeta_\sigma^n\cdot\mathbf{n},\zeta_c^n-\zeta_\lambda^n\rangle_{\partial\mathcal{T}_h}
+\langle\mathbf{u}_h^{n}\cdot\mathbf{n}\zeta_\lambda^n,\zeta_c^n-\zeta_\lambda^n\rangle_{\partial\mathcal{T}_h^{in}}+\langle\mathbf{u}_h^{n}\cdot\mathbf{n}\zeta_c^n,\zeta_c^n-\zeta_\lambda^n\rangle_{\partial\mathcal{T}_h^{out}}\\
=&(q\zeta_c^{n,\ast},\zeta_c^n)_{\mathcal{T}_h}+\frac{1}{\Delta t}(\phi\zeta_c^{n-1},\zeta_c^n)_{\mathcal{T}_h}.
\end{align*}

After rearranging and multiplying the above equation by $2\Delta t$, we can get  the following equation with \eqref{eq3-3}
\begin{equation}\nonumber
\begin{aligned}
&2(\phi(\zeta_c^n-\zeta_c^{n-1}),\zeta_c^n)_{\mathcal{T}_h}+2\Delta t(\nabla\zeta_c^n,\nabla\zeta_c^n)_{\mathcal{T}_h}+2\Delta t\langle\zeta_\lambda^n-\zeta_c^n,\frac{1}{h}(\zeta_\lambda^n-\zeta_c^n)\rangle_{\partial\mathcal{T}_h}
+\Delta t\langle|\mathbf{u}_h^{n}\cdot\mathbf{n}|(\zeta_\lambda^n-\zeta_c^n),\zeta_\lambda^n-\zeta_c^n\rangle_{\partial\mathcal{T}_h}\\
=&2\Delta t(\zeta_\sigma^n,\nabla \zeta_c^n)_{\mathcal{T}_h}-2\Delta t(D(\mathbf{u}_h^{n})^{-1}\zeta_\sigma^n,\nabla \zeta_c^n)_{\mathcal{T}_h}-2\Delta t\langle\zeta_\sigma^n\cdot\mathbf{n},\zeta_c^n-\zeta_\lambda^n\rangle_{\partial\mathcal{T}_h}-\Delta t(q\zeta_c^n,\zeta_c^n)_{\mathcal{T}_h}+2\Delta t(q\zeta_c^{n,\ast},\zeta_c^n)_{\mathcal{T}_h}\\
:=&F_1+F_2+F_3+F_4+F_5.
\end{aligned}
\end{equation}

For the first term of the left-hand-side of the above equation, we can get by the identity $2a(a-b)=a^2-b^2+(a-b)^2$
\begin{equation}\nonumber
\begin{aligned}
&2(\phi(\zeta_c^n-\zeta_c^{n-1}),\zeta_c^n)_{\mathcal{T}_h}\geq(\phi\zeta_c^n,\zeta_c^n)_{\mathcal{T}_h}-(\phi\zeta_c^{n-1},\zeta_c^{n-1})_{\mathcal{T}_h}.\\
\end{aligned}
\end{equation}

By Cauchy and Young's inequalities, we can reach the bound of $F_1, F_2,\ldots,F_5$ one by one as follows
\begin{equation}\nonumber
\begin{aligned}
&F_1\leq C\Delta t(\epsilon\|\nabla \zeta_c^n\|_{\mathcal{T}_h}^2+\frac{1}{\epsilon}\|\zeta_\sigma^n\|_{\mathcal{T}_h}^2),\\
&F_2\leq C\Delta t(\epsilon\|\nabla \zeta_c^n\|_{\mathcal{T}_h}^2+\frac{1}{\epsilon}\|\zeta_\sigma^n\|_{\mathcal{T}_h}^2),\\
&F_3\leq C\Delta t(\frac{\epsilon}{h}|\zeta_\lambda^n-\zeta_c^n|_{\partial\mathcal{T}_h}^2+\frac{1}{\epsilon}\|\zeta_\sigma^n\|_{\mathcal{T}_h}^2),\\
&F_4+F_5\leq C\Delta t\|\zeta_c^n\|_{\mathcal{T}_h}^2.\\
\end{aligned}
\end{equation}
Combining the above estimates, for sufficiently small $\epsilon$, we can derive
\begin{equation}\nonumber
\begin{aligned}
\|\zeta_c^n\|_{\mathcal{T}_h}^2-\|\zeta_c^{n-1}\|_{\mathcal{T}_h}^2+\Delta t\|\nabla \zeta_c^n\|_{\mathcal{T}_h}^2+\frac{\Delta t}{h}|\zeta_\lambda^n-\zeta_c^n|_{\partial\mathcal{T}_h}^2+\Delta t|\zeta_\lambda^n-\zeta_c^n|_{u,\partial\mathcal{T}_h}^2
\leq C\Delta t(\|\zeta_\sigma^n\|_{\mathcal{T}_h}^2+\|\zeta_c^n\|_{\mathcal{T}_h}^2)
\end{aligned}
\end{equation}
where $|v|_{u,\partial\mathcal{T}_h}^2:=\langle|\mathbf{u}_h^{n}\cdot\mathbf{n}|v,v\rangle_{\partial\mathcal{T}_h}$.

Thus, summing the above equation from $n=1$ to $m$ and using discrete Gronwal's inequality with $\zeta_c^n=0$, we can obtain
\begin{equation}\nonumber
\begin{aligned}
\|\zeta_c^m\|_{\mathcal{T}_h}^2+\Delta t\sum_{n=1}^m(\|\nabla \zeta_c^n\|_{\mathcal{T}_h}^2+\frac{1}{h}|\zeta_\lambda^n-\zeta_c^n|_{\partial\mathcal{T}_h}^2+|\zeta_\lambda^n-\zeta_c^n|_{u,\partial\mathcal{T}_h}^2)\leq C\Delta t\sum_{n=1}^m\|\zeta_\sigma^n\|_{\mathcal{T}_h}^2.
\end{aligned}
\end{equation}

\end{proof}

\begin{lemma}\label{err_b}
If the regularity assumption \ref{eq4} holds, there exists the following error bound
 \begin{equation}\nonumber
\begin{aligned}
h|\lambda_h^{c,n}-c^{n}|_{\partial\mathcal{T}_h}^2\leq& C(h^2\|\mathbf{u}^{n}-\mathbf{u}^{n}_h\|_{\mathcal{T}_h}^2+h^2\|\bm\sigma^{n}-\bm\sigma^{n}_h\|_{\mathcal{T}_h}^2+\|c^{n}-c_h^{n}\|_{\mathcal{T}_h}^2+h|c^{n}-\Pi_ec^{n}|_{\partial\mathcal{T}_h}^2),
\end{aligned}
\end{equation}
where  $C>0$ denotes a constant independent of the parameters $h$ and $\Delta t$.
\end{lemma}
\begin{proof}
Set $\bm\tau_h\cdot\mathbf{n}=\lambda_h^{c,n}-\Pi_ec^{n}$. A simple argument shows that
\begin{equation}\nonumber
\begin{aligned}
h\|\bm\tau_h\|_{1,\mathcal{T}_h}+\|\bm\tau_h\|_{\mathcal{T}_h}\leq Ch^{\frac{1}{2}}|\lambda_h^{c,n}-\Pi_ec^{n}|_{\partial\mathcal{T}_h}.\\
\end{aligned}
\end{equation}

Green's formula and $\bm\sigma+D(\mathbf{u})\nabla c=0$ imply
\begin{equation}\nonumber
\begin{aligned}
&(D(\mathbf{u}^{n})^{-1}\bm\sigma^{n},\bm\tau_h)_{\mathcal{T}_h}-(c^{n},\nabla\cdot\bm\tau_h)_{\mathcal{T}_h}+\langle c^{n},\bm\tau_h\cdot\mathbf{n}\rangle_{\partial\mathcal{T}_h}=0\\
\end{aligned}
\end{equation}
and
\begin{equation}\nonumber
\begin{aligned}
&(D(\mathbf{u}_h^{n})^{-1}\bm\sigma^{n}_h,\bm\tau_h)_{\mathcal{T}_h}-(c^{n}_h,\nabla\cdot\bm\tau_h)_{\mathcal{T}_h}+\langle \lambda_h^{c,n},\bm\tau_h\cdot\mathbf{n}\rangle_{\partial\mathcal{T}_h}=0.
\end{aligned}
\end{equation}
Hence we can get
\begin{equation}\nonumber
\begin{aligned}
&|\lambda_h^{c,n}-\Pi_ec^{n}|_{\partial\mathcal{T}_h}^2=\langle\lambda_h^{c}-\Pi_ec,\lambda_h^{c}-c\rangle_{\partial\mathcal{T}_h}\\
=&(D(\mathbf{u}^{n})^{-1}\bm\sigma-D(\mathbf{u}^{n}_h)^{-1}\bm\sigma_h,\bm\tau_h)_{\mathcal{T}_h}-(\nabla\cdot\bm\tau_h,c^{n}-c_h^{n})_{\mathcal{T}_h}\\
=&(\bm\sigma^{n}(D(\mathbf{u}^{n})^{-1}-D(\mathbf{u}^{n}_h)^{-1}),\bm\tau_h)_{\mathcal{T}_h}+(D(\mathbf{u}_h^{n})^{-1}(\bm\sigma^{n}-\bm\sigma^{n}_h),\bm\tau_h)_{\mathcal{T}_h}\\
&-(\nabla\cdot\bm\tau_h,c^{n}-c_h^{n})_{\mathcal{T}_h}.
\end{aligned}
\end{equation}
Thus, we can obtain
\begin{equation}\nonumber
\begin{aligned}
h^{\frac{1}{2}}|\lambda_h^{c,n}-\Pi_ec^{n}|_{\partial\mathcal{T}_h}\leq&C(h\|\mathbf{u}^{n}-\mathbf{u}^{n}_h\|_{\mathcal{T}_h}+h\|\bm\sigma^{n}-\bm\sigma^{n}_h\|_{\mathcal{T}_h}+\|c^{n}-c_h^{n}\|_{\mathcal{T}_h}).
\end{aligned}
\end{equation}

Using the triangle inequality, we get
\begin{equation}\nonumber
\begin{aligned}
h^{\frac{1}{2}}|\lambda_h^{c,n}-c^{n}|_{\partial\mathcal{T}_h}\leq C(h\|\mathbf{u}^{n}-\mathbf{u}^{n}_h\|_{\mathcal{T}_h}+h\|\bm\sigma^{n}-\bm\sigma^{n}_h\|_{\mathcal{T}_h}+\|c^{n}-c_h^{n}\|_{\mathcal{T}_h}+h^{\frac{1}{2}}|c^{n}-\Pi_ec^{n}|_{\partial\mathcal{T}_h}).
\end{aligned}
\end{equation}
\end{proof}

Now, we can complete our proof of the convergence theorem.

\begin{proof}
Set
\begin{equation}\nonumber
\begin{aligned}
\eta_c=\Pi_h c-c,\quad\eta_\sigma=\Pi^{RT} \bm\sigma-\bm\sigma,\quad\eta_\mathbf{u}=\Pi^{RT} \mathbf{u}-\mathbf{u}.
\end{aligned}
\end{equation}
Using \eqref{eq7}, \eqref{con_pro_3} and \eqref{con_pro_4}, we can get the error residual equation
\begin{equation}\label{eq4-1}
\begin{aligned}
(\up{a})\quad&(D(\mathbf{u}^{n})^{-1}\bm\sigma^{n}-D(\mathbf{u}^{n}_h)^{-1}\bm\sigma_h^{n},\bm\tau_h)_{\mathcal{T}_h}-(c^{n}-c_h^{n},\nabla\cdot\bm\tau_h)_{\mathcal{T}_h}=\langle\lambda_h^{c,n},\bm\tau_h\cdot\mathbf{n}\rangle_{\partial\mathcal{T}_h},\\
(\up{b})\quad&(\phi\partial_tc^{n}-\phi\bar{\partial}c^{n}_{h},v_h)_{\mathcal{T}_h}+(\nabla\cdot(\mathbf{u}^{n}c^{n})-\nabla\cdot(\mathbf{u}^{n}_hc^{n}_h),v_h)_{\mathcal{T}_h}\\
&\quad+\langle\mathbf{u}^{n}\cdot\mathbf{n}(\lambda^{c,{n}}-c^{n})-\mathbf{u}_h^{n}\cdot\mathbf{n}(\lambda_h^{c,{n}}-c_h^{n}),v_h\rangle_{\partial\mathcal{T}_h^{in}}\\
&\quad\quad+(\nabla\cdot\bm\sigma^{n}-\nabla\cdot\bm\sigma_h^{n},v_h)_{\mathcal{T}_h}=(\phi\partial_tc^{n}-\phi\f{\p c}{\p t},v_h)_{\mathcal{T}_h}+(q(c^{*,n}-c^{*,n}_{h}),v_h)_{\mathcal{T}_h},
\end{aligned}
\end{equation}
where $\partial_tc^{n}=(c^{n}-c^{n-1})/\Delta t$.

Taking $\bm\tau_{h}=\zeta^{n}_\sigma$ and $v_h=\zeta^{n}_c$ as test functions in \eqref{eq4-1} and using \eqref{con_pro_1} with $\mu_h=\Pi_ec^{n}-\lambda_h^{c,n}$  , we can reach that
 \begin{equation}\nonumber
\begin{aligned}
&(D(\mathbf{u}^{n})^{-1}\bm\sigma^{n}-D(\mathbf{u}^{n}_h)^{-1}\bm\sigma^{n}_h,\zeta^{n}_\sigma)_{\mathcal{T}_h}+ (\phi\partial_tc^{n}-\phi\bar{\partial}c^{n},\zeta^{n}_c)_{\mathcal{T}_h}\\
&+\langle\mathbf{u}^{n}\cdot\mathbf{n}(\lambda^{c,{n}}-c^{n})-\mathbf{u}_h^{n}\cdot\mathbf{n}(\lambda_h^{c,{n}}-c_h^{n}),\zeta^{n}_c\rangle_{\partial\mathcal{T}_h^{in}}+(\nabla\cdot(\mathbf{u}^{n}c^{n})-\nabla\cdot(\mathbf{u}_h^{n}c_h^{n}),\zeta^{n}_c)_{\mathcal{T}_h}
\\
=&(\phi\partial_tc^{n}-\phi\f{\p c}{\p t},\zeta^{n}_c)_{\mathcal{T}_h}+(q(c^{*,n}-c^{*,n}_{h}),\zeta^{n}_c)_{\mathcal{T}_h}.
\end{aligned}
\end{equation}
Multiplying the above equation by $2\Delta t$, we have
\begin{equation}\label{eq4-2}
\begin{aligned}
&2\Delta t((D(\mathbf{u}_h^{n})^{-1}\zeta^{n}_\sigma,\zeta^{n}_\sigma)_{\mathcal{T}_h}+2(\phi(\zeta^{n}_c-\zeta^{n-1}_c),\zeta^{n}_c)_{\mathcal{T}_h}\\
=&-2\Delta t((D(\mathbf{u}^{n})^{-1}-D(\mathbf{u}^{n}_h)^{-1})\bm\sigma^{n},\zeta^{n}_\sigma)_{\mathcal{T}_h}+2\Delta t(D(\mathbf{u}_h^{n})^{-1}\eta^{n}_\sigma,\zeta^{n}_\sigma)_{\mathcal{T}_h}\\
&+2(\phi(\eta^{n}_c-\eta^{n-1}_c),\zeta^{n}_c)_{\mathcal{T}_h}
+2\Delta t(\nabla\cdot(\mathbf{u}^{n}c^{n})-\nabla\cdot(\mathbf{u}^{n}_hc_h^{n}),\zeta^{n}_c)_{\mathcal{T}_h}\\
&+2\Delta t\langle\mathbf{u}^{n}\cdot\mathbf{n}(\lambda^{c,{n}}-c^{n})-\mathbf{u}_h^{n}\cdot\mathbf{n}(\lambda_h^{c,{n}}-c_h^{n}),\zeta^{n}_c\rangle_{\partial\mathcal{T}_h^{in}}
\\&+2\Delta t(\phi\partial_tc^{n}-\phi\f{\p c}{\p t},\zeta^{n}_c)_{\mathcal{T}_h}+2\Delta t(q^{n}(c^{*,n}-c^{*,n}_{h}),\zeta^{n}_c)_{\mathcal{T}_h}\\
=&Q_1+Q_2+Q_3+Q_4+Q_5+Q_{6}+Q_{7}.
\end{aligned}
\end{equation}

For the terms on the left hand side  of \eqref{eq4-2}, by the Lemma \ref{DU} and the identity $2a(a-b)=a^2-b^2+(a-b)^2$, we find
\begin{equation}\label{eq4-3}
\begin{aligned}
&2\Delta t((D(\mathbf{u}_h^{n})^{-1}(\zeta^{n}_\sigma),\zeta^{n}_\sigma)_{\mathcal{T}_h}\geq C\Delta t\|\zeta^{n}_\sigma\|_{\mathcal{T}_h}^2,\\
&2(\phi(\zeta^{n}_c-\zeta^{n-1}_c),\zeta^{n}_c)_{\mathcal{T}_h}\geq(\phi\zeta^{n}_c,\zeta^{n}_c)_{\mathcal{T}_h}-(\phi\zeta^{n-1}_c,\zeta^{n-1}_c)_{\mathcal{T}_h}.
\end{aligned}
\end{equation}

Next, we estimate the bound of $Q_1$, $Q_2$, $\ldots,$ $Q_7$ one by one. We first bound the terms $Q_1$, $Q_2$, $Q_3$, $Q_6$ and $Q_7$. By using Cauchy inequality, Young's inequality and Lemma \ref{err_u}, we have the following results
\begin{equation}\label{eq4-4}
\begin{aligned}
&Q_1\leq \Delta t\epsilon\|\zeta^{n}_\sigma\|_{\mathcal{T}_h}^2+C\Delta t\frac{1}{\epsilon}(\|\zeta^{n}_c\|_{\mathcal{T}_h}^2+\|\eta^{n}_c\|_{\mathcal{T}_h}^2+\|\eta^{n}_\mathbf{u}\|_{\mathcal{T}_h}^2),\\
&Q_2\leq \Delta t\epsilon\|\zeta^{n}_\sigma\|_{\mathcal{T}_h}^2+C\Delta t\frac{1}{\epsilon}\|\eta_\sigma^{n}\|_{\mathcal{T}_h}^2,\\
&Q_3\leq C(\Delta t\|\zeta^{n}_c\|_{\mathcal{T}_h}^2+\int^{t^{n}}_{t^{n-1}}\|\f{\p\eta_c}{\p t}\|_{\mathcal{T}_h}^2dt),\\
&Q_{6}\leq C(\Delta t^{2}\int^{t^{n}}_{t^{n-1}}\|\f{\p^{2}c}{\p t^{2}}\|_{\mathcal{T}_h}^2dt+\Delta t\|\zeta^{n}_c\|_{\mathcal{T}_h}^2),\\
&Q_7\leq C\Delta t(\|\zeta^{n}_c\|_{\mathcal{T}_h}^2+\|\eta^{n}_c\|_{\mathcal{T}_h}^2).
\end{aligned}
\end{equation}

And then, we show the boundedness of the terms $Q_4$ and $Q_5$. Taking $\mu_h=\zeta^{n}_c$ in \eqref{con_pro_2}, we can get
\begin{equation}\nonumber
\begin{aligned}
\langle\mathbf{u}_h^{n}\cdot\mathbf{n}\hat{c}_h^{n},\zeta^{n}_c\rangle_{\partial\mathcal{T}_h}=-\langle|\mathbf{u}_h^{n}\cdot\mathbf{n}|\lambda_h^{c,{n}},\zeta^{n}_c\rangle_{\partial\mathcal{T}_h^{in}}+\langle|\mathbf{u}_h^{n}\cdot\mathbf{n}|c_h^{n},\zeta^{n}_c\rangle_{\partial\mathcal{T}_h^{out}}=0.
\end{aligned}
\end{equation}
Thus we can reach
\begin{equation}\nonumber
\begin{aligned}
-\langle\mathbf{u}_h^{n}\cdot\mathbf{n}(\lambda_h^{c,{n}}-c_h^{n}),\zeta^{n}_c\rangle_{\partial\mathcal{T}_h^{in}}
=&\langle|\mathbf{u}_h^{n}\cdot\mathbf{n}|(\lambda_h^{c,{n}}-c_h^{n}),\zeta^{n}_c\rangle_{\partial\mathcal{T}_h^{in}}\\
=&\langle|\mathbf{u}_h^{n}\cdot\mathbf{n}|c_h^{n},\zeta^{n}_c\rangle_{\partial\mathcal{T}_h^{out}}-\langle|\mathbf{u}_h^{n}\cdot\mathbf{n}|c_h^{n},\zeta^{n}_c\rangle_{\partial\mathcal{T}_h^{in}}\\
=&\langle\mathbf{u}_h^{n}\cdot\mathbf{n}c_h^{n},\zeta^{n}_c\rangle_{\partial\mathcal{T}_h^{out}}+\langle\mathbf{u}_h^{n}\cdot\mathbf{n}c_h^{n},\zeta^{n}_c\rangle_{\partial\mathcal{T}_h^{in}}
=\langle\mathbf{u}_h^{n}\cdot\mathbf{n}c_h^{n},\zeta^{n}_c\rangle_{\partial\mathcal{T}_h}.
\end{aligned}
\end{equation}

Now we can estimate the terms $Q_4$ and $Q_5$ by Green's formula and Lemma \ref{err_u}
\begin{equation}\label{eq4-5}
\begin{aligned}
&Q_4+Q_5\\
=&2\Delta t(\nabla\cdot(\mathbf{u}^{n}c^{n})-\nabla\cdot(\mathbf{u}_h^{n}c_h^{n}),\zeta^{n}_c)_{\mathcal{T}_h}-2\Delta t\langle\mathbf{u}^n\cdot\mathbf{n}-\mathbf{u}^n_h\cdot\mathbf{n},\zeta^{n}_c\rangle_{\partial\mathcal{T}_h}\\
=&-2\Delta t(\mathbf{u}^{n}c^{n}-\mathbf{u}_h^{n}c_h^{n},\nabla\zeta^{n}_c)_{\mathcal{T}_h}
=-2\Delta t(\mathbf{u}^{n}_{h}(\zeta^{n}_c-\eta^{n}_c)+(\mathbf{u}^n-\mathbf{u}^n_h)c^{n},\nabla\zeta^{n}_c)_{\mathcal{T}_h}\\
=&-2\Delta t(\mathbf{u}^{n}(\zeta^{n}_c-\eta^{n}_c)+(\mathbf{u}^n-\mathbf{u}^n_h)c^{n},\nabla\zeta^{n}_c)_{\mathcal{T}_h}
+2\Delta t([\mathbf{u}^{n}-\mathbf{u}^{n}_h](\zeta^{n}_c-\eta^{n}_c)+(\mathbf{u}^n-\mathbf{u}^n_h)c^{n},\nabla\zeta^{n}_c)_{\mathcal{T}_h}\\
\leq&C\Delta t([1+h^{-d}\|\zeta^{n-1}_c\|_{\mathcal{T}_h}^2]\|\zeta^{n}_c\|_{\mathcal{T}_h}^2+\|\eta^{n}_c\|_{\mathcal{T}_h}^2
+\|\zeta^{n-1}_c\|_{\mathcal{T}_h}^2+\|\eta^{n}_{\mathbf{u}}\|_{\mathcal{T}_h}^2)+\Delta t\epsilon\|\nabla\zeta^{n}_c\|_{\mathcal{T}_h}^2.
\end{aligned}
\end{equation}
where we have used the fact that $\|\mathbf{u}^{n}\|_{L^{\infty}}, \|c^{n}\|_{L^{\infty}}\leq C$.

Finally, by collecting the estimates and using Lemma \ref{L2_err} and \ref{rt}, we can obtain
\begin{equation}\label{eq4-6}
\begin{aligned}
\|\zeta^{n}_c\|_{\mathcal{T}_h}^2-\|\zeta^{n-1}_c\|_{\mathcal{T}_h}^2+\Delta t\|\zeta^{n}_\sigma\|_{\mathcal{T}_h}^2\leq C\Delta t(h^{2s}+\Delta t^2+[1+h^{-d}\|\zeta^{n-1}_c\|_{\mathcal{T}_h}^2]\|\zeta^{n}_c\|_{\mathcal{T}_h}^2+\|\zeta^{n-1}_c\|_{\mathcal{T}_h}^2+\epsilon\|\nabla\zeta^{n}_c\|_{\mathcal{T}_h}^2+\epsilon\|\zeta^{n}_\sigma\|_{\mathcal{T}_h}^2),
\end{aligned}
\end{equation}

In order to complete our proof, we need the following inductive hypothesis
\begin{equation}\label{hypo}
h^{-d/2}\|\zeta^{n}_c\|_{\mathcal{T}_h}\leq C,\quad  \forall n\geq 0.
\end{equation}
Assume that the inductive hypothesis \eqref{hypo} holds for $n=0,1,\cdots,m-1$. When $n=m$,  summing \eqref{eq4-6} from $1$ to $m$, for  sufficiently small  $\epsilon$, using the discrete Gronwall's inequality with $\xi^{0}_c=0$ and Lemma \ref{err_gradc}, we can get 
\begin{equation}\label{eq31}
\|\zeta^{m}_c\|_{\mathcal{T}_h}^2+\Delta t\sum\limits_{n=0}^m\|\zeta^{n}_\sigma\|_{\mathcal{T}_h}^2\leq C(h^{2s}+\Delta t^2).
\end{equation}
Thus, the estimate \eqref{eq24}(d) follows from \eqref{eq24}(b)  and Lemma \ref{err_b}. As in \cite{model2}, we know that
\[
\|\Pi_hp^{n}-p_h^{n}\|_{\mathcal{T}_h}\leq C(1+\|\Pi^{RT}\mathbf{u}^{n}\|_{L^{\infty}})\|c_h^{n}-c^{n}\|_{\mathcal{T}_h},
\]
so  the estimates \eqref{eq24}(a) and \eqref{eq24}(c) can be straight obtained by Lemma \ref{err_u}.

As we know, the error estimate \eqref{eq31} is obtained under the inductive hypothesis \eqref{hypo}. Now we check it. When $n=0$, we know that $\zeta^0_c=0$. So the inductive hypothesis \eqref{hypo} holds. When $n=m$, from \eqref{eq31} we know that
\begin{equation}\nonumber
h^{-d/2}\|\zeta^{m}_c\|_{\mathcal{T}_h}\leq Ch^{-d/2}(h^{s}+\Delta t)\leq C.
\end{equation}
Hence the induction hypothesis \eqref{hypo} holds for $n=m$.
\end{proof}

\section*{Acknowledgments}
Zhang's work was supported partially by  the Natural Science Foundation of Shandong Province (ZR2019MA015) and the Fundamental Research Funds for the Central Universities (22CX03020A). Zhu’s work was partially supported by the National Council for Scientific and Technological Development of Brazil (CNPq).


\end{document}